\documentclass[10pt]{amsart}
\usepackage{latexsym}
\usepackage{amsmath}
\usepackage{amssymb}
\usepackage{mathrsfs}
\usepackage{graphicx}
\usepackage{bbm}
\usepackage{bbold}
\input epsf
\input xy
\xyoption{all}

\def\Orth{\mathop{\fam0 Orth}\nolimits}
\def\dom{\mathop{\fam0 dom}\nolimits}
\def\On{\mathop{\fam0 On}\nolimits}

\begin{document}
\title[Boolean Models and Simultaneous Inequalities ]
{Boolean Models and\\ Simultaneous Inequalities}
\author{S. S. Kutateladze}
\date{August 4, 2009}
\address[]{
Sobolev Institute of Mathematics\newline
\indent 4 Koptyug Avenue\newline
\indent Novosibirsk, 630090\newline
\indent RUSSIA
}
\begin{abstract}
Boolean valued analysis is applied to deriving  operator
versions of the classical Farkas Lemma in the  theory of simultaneous
linear inequalities.
\end{abstract}
\email{
sskut@member.ams.org
}
\thanks{This talk was delivered at the Maltsev Centennial Meeting
(Novosibirsk, Russia) on August 27, 2009.}

\maketitle

\section{Introduction}
The Farkas Lemma, also known as the Farkas--Minkowski Lemma,
plays a key role in linear programming and the relevant areas of
optimization.\footnote{Cp.~\cite{History}\cite{Encyc}.}
The aim of this talk is to demonstrate how Boolean valued
analysis\footnote{Cp.~\cite{IBA}.}  may
be applied to simultaneous linear inequalities with operators.
This particular theme is another illustration of the deep and powerful
technique of ``stratified validity'' which is characteristic of
Boolean valued analysis.

Assume that $X$ is a~real vector space,
$Y$ is a~{\it Kantorovich space\/} also known as a~Dedekind complete vector lattice or a complete Riesz
space. Let $\mathbb  B:=\mathbb B(Y)$ be the {\it base\/} of~$Y$, i.e., the complete Boolean algebras of
positive projections in~$Y$; and let $m(Y)$ be the  universal completion of~$Y$.
Denote by $L (X,Y)$ the space of linear operators from $X$ to~$Y$
and let $\Orth(Y)$ stand for the commutant of $\mathbb  B$ 
in $L^{(r)}(Y)$. 
In case $X$ is furnished with some $Y$-seminorm on~$X$, by $L^{(m)}(X,Y)$ we mean
the {\it space of dominated operators\/} from $X$ to~$Y$.
As usual,  $\{T\le y\}:=\{x\in X \mid Tx\le y\}$ and $\ker(T)=T^{-1}(0)$
for $T: X\to Y$.

\section{Inequalities: Explicit Dominance}

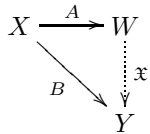
\[
\xymatrix{
  X\ar[dr]_{B} \ar[r]^{A}
                & W  \ar@{.>}[d]^{\mathfrak X}  \\
                & Y           }
\]

{(\bf 1):}
$(\exists \mathfrak X)\ {\mathfrak X}A=B \leftrightarrow {\ker(A)\subset\ker(B)};
$

{(\bf 2):}\footnote{The Kantorovich  Theorem, \cite[p.~44]{SubDif}.} {\sl If $W$ is ordered by $W_+$ and $A(X)-W_+=W_+ - A(X)=W$, then}
$$
(\exists \mathfrak X\ge 0)\ {\mathfrak X}A=B \leftrightarrow \{A\le0\}\subset\{B\le 0\}.
$$

\section{Farkas: Explicit Dominance}
{\scshape Theorem~1.}{\sl\
Assume  that $A_1,\dots,A_N$ and $B$ belong to~$ L^{(m)}(X,Y)$.

The following are equivalent:

{\rm(1)} Given   $b\in \mathbb B$, the operator inequality $bBx\le 0$
is a consequence of the simultaneous linear operator inequalities
$bA_1 x\le0,\dots, bA_N x\le 0$,
i.e.,
$$
\{bB\le 0\}\supset\{bA_1\le 0\}\cap\dots\cap\{bA_N\le 0\}.
$$

{\rm(2)} There are  positive  orthomorphisms $\alpha_1,\dots,\alpha_N\in\Orth(m(Y))$
such that
$$
B=\sum\limits_{k=1}^N{\alpha_k A_k};
$$
 i.e.,
$B$ lies in the operator convex
conic hull of~$A_1,\dots,A_N$.}

\section{Farkas: Hidden Dominance}
{\scshape Lemma 1.}{\sl\
Let $X$ be a vector space over some subfield $R$ of the reals~
$\mathbb{R}$.
Assume that $f$ and $g$ are  $R$-linear functionals on~$X$;
in symbols, $f, g\in X^{\#}:=L(X,\mathbb R)$.

For the inclusion
$$
\{g\le0\}\supset\{f\leq0\}
$$
to hold it is necessary and sufficient that
there be $\alpha\in\mathbb{R}_+$ satisfying
$g=\alpha f$.}

{\scshape Proof.}
{\scshape Sufficiency} is obvious.

{\scshape Necessity:} The case of $f=0$ is trivial. If $f\ne0$ then there is
some $x\in X$ such that $f(x)\in \mathbb R$ and~$f(x)>0$. Denote the image $f(X)$ of~$X$
under~$f$ by~$R_0$. Put $h:=g\circ f^{-1}$, i~.e. $h\in R_0^{\#}$ is the only solution
for $h\circ f=g$. By hypothesis, $h$ is a positive $R$-linear functional
on~$R_0$. By the Bigard Theorem \cite[p.~108]{SubDif}
$h$ can be extended to a positive homomorphism
$\bar{h}: \mathbb R\to\mathbb R$, since
$R_0-\mathbb R_+=\mathbb R_+-R_0=\mathbb R$.
Each positive automorphism of~$\mathbb R$ is multiplication by a positive real.
As the sought $\alpha$ we may take  $\bar{h}(1)$.

The proof of the lemma is complete.

\section{ Reals: Explicit Dominance}
{\scshape Lemma 2.}{\sl\
Let $X$ be an $\mathbb R$-seminormed vector space over some subfield $R$ of~
$\mathbb{R}$.
Assume that $f_1,\dots, f_N$ and $g$ are  bounded $R$-linear functionals on~$X$; in symbols, $f_1,\dots, f_N, g\in X^{*}:=L^{(m)}(X,\mathbb R)$.

For the inclusion
$$
\{g\le0\}\supset\bigcap\limits_{k=1}^{N}\{f_k\leq0\}
$$
to hold it is necessary and sufficient that
there be $\alpha_1,\dots,\alpha_N \in\mathbb{R}_+$ satisfying
$$
g=\sum\limits_{k=1}^{N}\alpha_k f_k.
$$
}

\section{Origins}
Cohen's final solution of the problem of the cardinality of the
continuum within ZFC gave rise to the Boolean valued models.

Scott forecasted  in 1969:\footnote{Cp.~\cite{Scott}.}
\begin{itemize}
\item[]{{\it We must ask whether there is any interest
in these nonstandard models aside from the independence proof;
that is, do they have any mathematical interest?
The answer must be yes, but we cannot yet give a really good argument.}
}
\end{itemize}

Takeuti coined the term ``Boolean valued
analysis''  for applications of the models
to analysis.\footnote{Cp.~\cite{Takeuti}.}

\section{Boolean Valued Universe}
Let
${\mathbb B}$
be a~complete Boolean algebra. Given an ordinal
$\alpha$,
put
$$
%\gathered
V_{\alpha}^{({\mathbb B})}
:=\{x \mid
%\\
(\exists \beta\in\alpha)\ x:\dom (x)\rightarrow
{\mathbb B}\ \&\ \dom (x)\subset V_{\beta}^{({\mathbb B})}  \}.
%\endgathered
$$
The {\it Boolean valued
universe\/}
${\mathbb V}^{({\mathbb B})}$
is
$$
{\mathbb V}^{({\mathbb B})}:=\bigcup\limits_{\alpha\in\On} V_{\alpha}^{({\mathbb B})},
$$
with $\On$ the class of all ordinals. The truth
value $[\![\varphi]\!]\in {\mathbb B}$ is assigned to each formula
$\varphi$ of ZFC relativized to ${\mathbb V}^{({\mathbb B})}$.

\section{~Descending and Ascending}

Given $\varphi$, a~formula of ZFC, and
$y$, a~member of~${\mathbb V}^{{\mathbb B}}$; put
$A_{\varphi}:=A_{\varphi(\cdot,\ y)}:=\{x\mid\varphi (x,\ y)\}$.
The {\it descent\/}
$A_{\varphi}{\downarrow}$
of a~class
$A_{\varphi}$
is
$$
A_{\varphi}{\downarrow}:=\{t\mid t\in {\mathbb V}^{({\mathbb B})} \ \&\  [\![\varphi  (t,\ y)]\!]=\mathbb 1\}.
$$
If
$t\in A_{\varphi}{\downarrow}$,
then
it is said
that
{\it $t$
satisfies
$\varphi (\cdot,\ y)$
inside
${\mathbb V}^{({\mathbb B})}$}.
The {\it descent\/}
$x{\downarrow}$
of $x\in {\mathbb V}^{({\mathbb B})}$
is defined as
$$
x{\downarrow}:=\{t\mid t\in {\mathbb V}^{({\mathbb B})}\ \&\  [\![t\in x]\!]=\mathbb 1\},
$$
i.e. $x{\downarrow}=A_{\cdot\in x}{\downarrow}$.
The class $x{\downarrow}$ is a~set.
If $x$ is a~nonempty set
inside
${\mathbb V}^{({\mathbb B})}$ then
$$
(\exists z\in x{\downarrow})[\![(\exists t\in x)\ \varphi (t)]\!]
=[\![\varphi(z)]\!].
$$
The {\it ascent\/} functor acts in the opposite direction.

\section{ The Reals Within}

There is an object
$\mathscr R$
inside
${\mathbb V}^{({\mathbb B})}$ modeling $\mathbb R$, i.e.,
$$
[\![\mathscr R\ {\text{is the reals}}\,]\!]=\mathbb 1.
$$
Let $\mathscr R{\downarrow}$ be the descent of
 the carrier $|\mathscr R|$ of the algebraic system
$$
\mathscr R:=(|\mathscr R|,+,\,\cdot\,,0,1,\le)
$$
inside ${\mathbb V}^{({\mathbb B})}$.
Implement the descent of the structures on $|\mathscr R|$
to $\mathscr R{\downarrow}$ as follows:
$$
\gathered
x+y=z\leftrightarrow [\![x+y=z]\!]=\mathbb 1;
\\
xy=z\leftrightarrow [\![xy=z]\!]=\mathbb 1;
\\
x\le y\leftrightarrow [\![x\le y]\!]=\mathbb 1;
\\
\lambda x=y\leftrightarrow [\![\lambda^\wedge x=y]\!]=\mathbb 1
\\
(x,y,z \in\mathscr R{\downarrow},\ \lambda\in\mathbb R).
\endgathered
$$

{\scshape Gordon Theorem.}
{\sl $\mathscr R{\downarrow}$
with the descended structures is a universally complete
vector lattice with base
$\mathbb B(\mathscr R{\downarrow})$
isomorphic to
${\mathbb B}$}.

{ Proof of Theorem~1.}

{(2)$\to$ (1):}
If  $B=\sum\nolimits_{k=1}^N{\alpha_k A_k}$ for some positive
$\alpha_1,\dots,\alpha_N$ in~$\Orth(m(Y))$ while  $bA_k x\le 0$ for $b\in\mathbb B$ and
$x\in X$, then
$$
bBx=b\sum\limits_{k=1}^N \alpha_k A_k x=\sum\limits_{k=1}^N \alpha_k bA_k x\le 0
$$
since orthomorphisms commute and projections are orthomorphisms of~$m(Y)$.

{(1)$\to$ (2):}
Consider the separated Boolean valued universe $\mathbb{V}^{(\mathbb{B})}$ over the base~$\mathbb B$
of~$Y$.
By the Gordon Theorem the ascent $Y{\uparrow}$ of~$Y$ is $\mathscr{R}$,
the  reals inside~$\mathbb{V}^{(\mathbb{B})}$.

Using the canonical embedding, we see that $X^{\scriptscriptstyle\wedge}$
is an~$\mathscr R$-seminormed vector space over the standard name  $\mathbb{R}^{\scriptscriptstyle\wedge}$
of the reals~$\mathbb{R}$.
Moreover, $\mathbb{R}^{\scriptscriptstyle\wedge}$ is a~subfield and
sublattice of~${\mathscr R}=Y{\uparrow}$
inside~$\mathbb{V}^{(\mathbb{B})}$.

Put $f_k:=A_k{\uparrow}$ for all $k:=1,\dots, N$  and $g:=B{\uparrow}$.
Clearly,  all $f_1,\dots, f_N,g$ belong to~$(X^{\scriptscriptstyle\wedge})^*$
 inside~${\mathbb V}^{\mathbb B}$.

Define the finite sequence
$$
f:
\{1,\dots,N\}^{\scriptscriptstyle\wedge}\to (X^{\scriptscriptstyle\wedge})^*
$$
  as the
ascent of~$(f_1,\dots,f_N)$.
In other words, the truth values are as follows:
$$
[\![f_{k^{\scriptscriptstyle\wedge}}(x^{\scriptscriptstyle\wedge})=A_k x]\!]={\mathbb 1},\quad [\![g(x^{\scriptscriptstyle\wedge})=Bx]\!]={\mathbb 1}
$$
for all $x\in X$ and $k:=1,\dots,N$.

Put
$$
b:=[\![A_1x\le 0^{\scriptscriptstyle\wedge}]\!]\wedge\dots\wedge
[\![A_N x\le 0^{\scriptscriptstyle\wedge}]\!].
$$
Then $bA_k x\le 0$ for all $k:=1,\dots,N$ and
$bBx\le 0$ by~(1).

 Therefore,
$$
[\![A_1x\le 0^{\scriptscriptstyle\wedge}]\!]\wedge\dots\wedge
[\![A_N x\le 0^{\scriptscriptstyle\wedge}]\!]\le [\![Bx\le 0^{\scriptscriptstyle\wedge}]\!].
$$

In other words,
$$
%\gathered
[\![(\forall k:=1^{\scriptscriptstyle\wedge},\dots,
N^{\scriptscriptstyle\wedge})f_{k}(x^{\scriptscriptstyle\wedge})
\le0^{\scriptscriptstyle\wedge}]\!]
%\\
$$
$$
=
\bigwedge_{k:=1,\dots,N} [\![f_{k^{\scriptscriptstyle\wedge}}(x^{\scriptscriptstyle\wedge})\le0^{\scriptscriptstyle\wedge}]\!]
\le [\![g(x^{\scriptscriptstyle\wedge})\le 0^{\scriptscriptstyle\wedge}]\!].
%\endgathered
$$

Using  Lemma~2 inside~$\mathbb{V}^{(\mathbb{B})}$  and  appealing  to the maximum principle
of Boolean valued analysis, we infer that there is
a~finite sequence  $\alpha:\{1^{\scriptscriptstyle\wedge},
\dots,N^{\scriptscriptstyle\wedge}\}
\to \mathscr R_+$ inside ${\mathbb V}^{({\mathbb B})}$
satisfying
$$
[\![(\forall x\in X^{\scriptscriptstyle\wedge})
g(x)
=\sum_{k=1^{\scriptscriptstyle\wedge}}^{N^{\scriptscriptstyle\wedge}}
\alpha(k)  f_k(x)]\!]={\mathbb1}.
$$
Put
$\alpha_k:=\alpha(k^{\scriptscriptstyle\wedge})
\in {\mathscr R}_+{\downarrow}$
for $k:=1,\dots,N$.
Multiplication by an element in~$\mathscr R{\downarrow}$ is an orthomorphism of~$m(Y)$.
Moreover,
$$
B=\sum\limits_{k=1}^N \alpha_k  A_k,
$$
 which completes the proof.

\section{ Counterexample: No Dominance }
Lemma~1, describing the consequences of a single inequality,
does not  restrict the class of functionals under consideration.
The analogous version of the Farkas Lemma simply fails for
two simultaneous inequalities in general. Indeed, the inclusion
$\{f=0\}\subset\{g\leq 0\}$
equivalent to the inclusion $\{f=0\}\subset\{g=0\}$ does not imply
that $f$ and~$g$ are proportional in the case of an arbitrary subfield
of~$\mathbb R$. It suffices to look at $\mathbb R$ over the rationals $\mathbb Q$,
take some  discontinuous $\mathbb Q$-linear functional on~$\mathbb Q$  and
the identity automorphism of~$\mathbb Q$.
This gives grounds for the next result.

\section{ Reconstruction: No Dominance }
{\scshape Theorem 2.}{\sl\
Take $A$ and  $B$ in $L(X,Y)$.
The following are equivalent:

{\rm(1)} $(\exists \alpha\in \Orth(m(Y)))\ B=\alpha A$;

{\rm(2)} There is a projection $\varkappa\in \mathbb B$ such that
$$
%\gathered
\{\varkappa bB\le 0\}\supset\{\varkappa bA\le 0\};\quad
%\\
 \{\neg\varkappa bB\le 0\}\supset\{\neg\varkappa bA\ge 0\}
%\endgathered
$$
for all $b\in\mathbb B$.}

{\scshape Proof.}
Boolean valued analysis reduces the claim to
the scalar case. Applying Lemma~1 twice and writing  down the truth values,
complete the proof.

\section{ Interval Operators}
Let $X$ be a~vector lattice.
An~{\it interval operator\/} ${\bf T}$ from $X$ to~$Y$
is  an~order interval $[\underline{T}, \overline{T}]$
in  $L^{(r)}(X,Y)$, with $\underline{T}\le \overline{T}$.

The interval equation ${\bf B}=\mathfrak X {\bf A}$ has a~{\it weak interval solution\/}\footnote{Cp.~\cite{Inexact}.}
provided that $(\exists \mathfrak X )(\exists A\in {\bf A})(\exists B\in {\bf B})\ B=\mathfrak X A$.

Given an interval operator ${\bf T}$ and $x\in X$, put
 $$
 P_{{\bf T}}(x)=\overline{T}x_+ -\underline{T}x_-.
 $$
Call ${\bf T}$ is {\it adapted\/} in case
$\overline{T}-\underline{T}$ is the sum of finitely many disjoint addends,
and put  $\sim(x):=-x$ for all $x\in X$.

\section{ Interval Equations}

{\scshape Theorem 3.}
{\sl Let $X$ be a vector lattice, and let $Y$ be a~Kantorovich space.
Assume that
 ${\bf A}_1,\dots, {\bf A}_N$ are adapted interval operators and ${\bf B}$ is an arbitrary
 interval operator in the space of order bounded operators $L^{(r)}(X,Y)$.

The following are equivalent:

{\rm(1)} The interval equation
$$
{\bf B}=\sum_{k=1}^N \alpha_k {\bf A}_k
$$
has a weak interval solution $\alpha_1,\dots,\alpha_N\in\Orth(Y)_+$.

{\rm(2)} For all  $b\in \mathbb B$ we have
$$
\{b\mathfrak B\ge 0\}\supset\{b{\mathfrak A}_1^{\sim}\le 0\}\cap\dots\cap\{b{\mathfrak A}_N^{\sim}\le 0\},
$$
where ${\mathfrak A}_k^{\sim}:=P_{{\bf A}_k}\circ\sim$ for $k:=1,\dots,N$ and $\mathfrak B:=P_{{\bf B}}$.
}

\section{ Inhomogeneous Inequalities}

{\scshape Theorem 4.}
{\sl Let $X$ be a $Y$-seminormed space, with $Y$ a~Kantorovich space.
Assume given some dominated operators
$A_1,\dots,A_N,  B\in L^{(m)}(X,Y)$ and elements $u_1,\dots, u_N,v\in Y$.
The following are equivalent:

{\rm(1)} For all   $b\in \mathbb B$ the  inhomogeneous operator inequality $bBx\le bv$
is a consequence of the
 consistent simultaneous inhomogeneous operator inequalities
 $bA_1 x\le bu_1,\dots, bA_N x\le bu_N$,
i.e.,
$$
\{bB\le bv\}\supset\{bA_1\le bu_1\}\cap\dots\cap\{bA_N\le bu_N\}.
$$

{\rm(2)} There are positive orthomorphisms $\alpha_1,\dots,\alpha_N\in\Orth(m(Y))$
satisfying
$$
B=\sum\limits_{k=1}^N{\alpha_k A_k};\quad v\ge \sum\limits_{k=1}^N{\alpha_k u_k}.
$$
}

\section{Inhomogeneous Matrix Inequalities}
In applications we encounter inhomogeneous matrix inequalities
over various finite-dimensional spaces.\footnote{Cp.~\cite{Olvi}.}

{\scshape Theorem 5.}{\sl\
Let $X$ be a~real $Y$-seminormed space, with $Y$  a~Kantorovich space.
Assume that
$A\in L^{(m)}(X,Y^s)$, $A\in L^{(m)}(X,Y^s)$, $u\in Y^t$ and $v\in Y^m$, where
$s$ and $t$ are some naturals.

The following are equivalent:

{\rm(1)} For all   $b\in \mathbb B$  the inhomogeneous operator inequality
$bBx\le bv$ is a~consequemce of the consistent inhomogeneous inequality
$bA x\le bu$, i.e., $\{bB\le bv\}\supset\{bA\le bu\}$.

{\rm(2)} There is some $s\times t$ matrix with entries positive orthomorphisms
of~$m(Y)$ such that
$B={\mathfrak X} A$ and~${\mathfrak X}u\le v$
for the corresponding linear operator $\mathfrak X\in L_+(Y^s, Y^t)$.
}

\section{ Complex Scalars}

{\scshape Theorem 6.}{\sl\
Let $X$ be a~complex $Y$-seminormed space, with $Y$ a~Kantorovich space.
Assume given  $u_1,\dots, u_N,v\in Y$ and dominated operators
$A_1,\dots,A_N,  B\in L^{(m)}(X,Y_{\mathbb C})$ from~$X$ into the complexification
$Y_{\mathbb C}:=Y\otimes iY$ of~$Y$.
The following are equivalent:

{\rm(1)} For all   $b\in \mathbb B$ and $x\in X$ the inhomogeneous  inequality
$b\vert Bx\vert\le bv$
is a~consequence of the consistent simultaneous inhomogeneous  inequalities
$b\vert A_1 x\vert\le bu_1,\dots,b\vert A_N x\vert\le bu_N$, i.e.,
$$
\{b\vert B\vert\le bv\}\supset\{b\vert A\vert_1\le bu_1\}\cap\dots\cap\{b\vert A\vert_N\le bu_N\}.
$$

{\rm(2)} There are complex orthomorphisms $c_1,\dots, c_N\in\Orth(m(Y)_{\mathbb C})$
satisfying
$$
B=\sum\limits_{k=1}^N{c_k A_k};\quad v\ge \sum\limits_{k=1}^N{\vert c_k\vert u_k}.
$$
}

\section{ Theorem of the Alternative}

{\scshape Theorem 7.}{\sl\
Let $X$ be a~$Y$-seminormed real vector space,
with $Y$ a~Kantorovich space.
Assume  that $A_1,\dots,A_N$ and $B$ belong to~$ L^{(m)}(X,Y)$.

Then one and  only one of the following holds:

{\rm(1)} There are   $x\in X$ and $b, b' \in \mathbb B$  such that
$b'\le b$ and
$$
 b'Bx>0, bA_1 x\le0,\dots, bA_N x\le 0.
$$

{\rm(2)} There are   $\alpha_1,\dots,\alpha_N\in\Orth(m(Y))_+$
such that
$
B=\sum\nolimits_{k=1}^N{\alpha_k A_k}.
$
}

\section{ All Is Number}
The above results, although curious to some extent,
are nothing more than simple illustrations of the powerful technique
of model theory shedding new light at the {\it Pythagorean
Thesis}.
The theory of the reals enriches  mathematics,
demonstrating the liberating role of logic.
%\end{document}

\section{Pursuit of Truth}
We definitely feel truth, but we cannot define truth properly.
That is what Tarski explained to us in the  1930s. 

We pursue truth by way of proof, as wittily phrased by Mac Lane. 
Model theory evaluates and counts truth and proof.

The chase of truth not only leads us close to the truth we pursue but also enables  
us to nearly catch up with many other 
instances of truth which we were not aware nor even foresaw at the start 
of the rally pursuit.  That is what we have learned from 
the Boolean models elaborated in the  1960s by  Scott, 
Solovay, and  Vop\v enka.

%\section{References}

\bibliographystyle{plain}

\end{document}